\documentclass[11pt]{mcom-l}
\usepackage{amssymb,amsfonts,amsmath,epsfig,hyperref}
\usepackage[margin=1.33in]{geometry}
\usepackage{mathrsfs}
\usepackage{color}     
\newtheorem{theorem}{Theorem}[section]
\newtheorem{lemma}[theorem]{Lemma}
\theoremstyle{remark}
\newtheorem{remark}{\bf Remark}[section]
\newcommand{\nn}{\nonumber}
\def\R{\mathbb{R}}

\def\u{{\bf u}}
\def\v{{\bf v}}
\def\w{{\bf w}}
\def\e{{\bf e}}
\def\d{\,{\rm d}}
\title{ 
A bounded numerical solution with\\ a small 
mesh size
implies existence of\\ a smooth solution 
to the Navier--Stokes equations 
} 

\author{Buyang Li}
\address{Department of Applied Mathematics, 
The Hong Kong Polytechnic University, Hung Hom, Hong Kong.} 
\email {\href{mailto:buyang.li@polyu.edu.hk}{buyang.li{\it @}polyu.edu.hk}}

\begin{document}

\maketitle

\begin{abstract}
\small  
We prove that for a given smooth initial value, if the 
finite element solution of the three-dimensional Navier--Stokes equations is bounded in a certain norm with a relatively small mesh size, then the solution of the Navier--Stokes equations with this given initial value must be smooth and unique, and is successfully approximated by the numerical solution.
\end{abstract}



\pagestyle{myheadings}
\markboth{}{}

\section{Introduction}
\setcounter{equation}{0}

The Navier--Stokes equations have wide applications in many scientific and engineering fields, such as ocean currents, weather forecast, and air flow around a wing. Regardless of the wide range of their applications, whether the Navier--Stokes equations always admit a unique smooth solution is not known yet in three-dimensional domains for general smooth initial data. Global existence of weak solutions to the Navier--Stokes equations was proved by Leray and Hopf \cite{Hopf51,Leray34}, while it was recently shown in \cite{Buckmaster-Vicol-2019} that weak solutions with finite kinetic energy are not unique in general. Many recent efforts have been made in proving global well-posedness for small initial data \cite{FK64,Kato84,KT01,LL11} and blowup examples for some related equations \cite{Cheskidov,GP09,HL09,MS01,Tao}.

Driven by the various applications, many numerical methods have been proposed for solving the Navier--Stokes equations, such as the finite element methods \cite{HR82,HR90,Rann00}, finite difference methods \cite{Chorin68}, spectral methods \cite{Heywood82,STW}, the Lagrange--Galerkin method \cite{BSS12,He-2013,Pironneau,Suli88}, and the projection method for time discretization \cite{Chorin,DGLQ82,Shen92}. The convergence of numerical solutions to the Navier--Stokes equations in three-dimensional domains was all proved by assuming that the equations have a sufficiently smooth solution. 
A natural question is, given a bounded numerical solution, what can we say about the smoothness of the true solution without assuming well-posedness of the Navier--Stokes equations?

We answer this question partially in this paper: for any given smooth initial value, if a numerical solution remains bounded in a certain norm for some mesh size which is smaller than some positive constant (determined by the bound of the numerical solution), then the Navier--Stokes equations have a unique smooth solution and, simultaneously, the numerical solution successfully approximates the true solution. Note that only one numerical solution is needed to draw the conclusion, instead of a sequence of numerical solutions with mesh size tending to zero. To illustrate the idea, we consider the Navier-Stokes equations
\begin{align}
& \partial_t\u+\u\cdot\nabla\u
-\mu\Delta \u+\nabla p=0 ,\label{NSPDE1}\\
&\nabla\cdot\u=0 , \label{NSPDE2}
\end{align}
in a convex polyhedron $\Omega\subset\R^3$ with the Dirichlet boundary condition $\u={\bf 0}$ on $\partial\Omega$ and a given initial condition $\u(x,0)=\u^0(x)$ (where $\mu>0$ is the viscosity constant), and focus on a simple linearized finite element method for the discretization of the Navier--Stokes equations. As usual, we impose the condition $\int_\Omega p(x,t)\d x=0$ for the uniqueness of pressure.

In this paper, we only provide a theoretical result and a basic framework to obtain such results. We hope that results that are useful in practical computation can be obtained in the future by different analysis.


\section{Notations and main results}
\setcounter{equation}{0}


For any nonnegative real number $k$, we denote by $H^k$ the conventional Sobolev space of functions defined on $\Omega$, with abbreviation $L^2=H^0$, and denote by $H^1_0$ the subspace of $H^1$ consisting of functions whose traces on the boundary are zero; see \cite{Adams}. 
We denote by $L^2_0$ the subspace of $L^2$ consisting of functions with vanishing integrals over $\Omega$. The following vector-valued spaces related to incompressible flow will be used in this paper:
\begin{align*}
&{\bf L}^p=(L^p)^3
\quad\mbox{and}\quad{\bf H}^k= (H^k)^3 , \\
&\dot {\bf D}
=\{ {\bf v} \in C^\infty_0(\Omega)^3 :
\nabla\cdot{\bf v}=0 \} , \\
%
&\dot {\bf L}^2 
= \mbox{The completion of $\dot {\bf D}$ in ${\bf L}^2$}, \\
&\dot {\bf H}^1 
= \mbox{The completion of $\dot {\bf D}$ in ${\bf H}^1$} ,\\
&
\dot {\bf H}^2 = \dot {\bf H}^1 \cap {\bf H}^2.
\end{align*}
For the simplicity of notation, we denote by $\|\cdot\|_{H^k}$ the norms of both $H^k$ and ${\bf H}^k$, and denote by $\|\cdot\|_{L^p}$ the norms of both $L^p$ and ${\bf L}^p$. We denote by $\dot{\bf H}^{-1}$ the dual space of $\dot {\bf H}^1$ and denote the norm of $\dot{\bf H}^{-1}$ by $\|\cdot\|_{\dot H^{-1}}$.

We denote by $(\cdot,\cdot)$ the inner product of both $L^2$ and ${\bf L^2}$, and denote by $P:{\bf L}^2\rightarrow\dot{\bf L}^2$ the $L^2$-orthogonal projection onto the divergence-free subspace. 
Let $D(A)=\dot {\bf H}^2 $ and let $A=-P\Delta: D(A)\rightarrow \dot {\bf L}^2$ be the operator defined by 
$$
(A\v,\w)=(\nabla \v,\nabla \w) \quad\forall\,\v\in\dot {\bf H}^2,\,\, \w\in \dot {\bf H}^1. 
$$
On a convex polyhedral domain $\Omega$ the regularity results of Stokes equations in \cite{Dauge-1989} imply that equation $A{\bf v}={\bf f}\in \dot{\bf L}^2$ has a unique solution ${\bf v}\in \dot{\bf H}^2$. 
Then, according to \cite[Theorem in \textsection 1.15.3]{Triebel}, the domain $D(A^{\frac{s}{2}})$ coincides with the complex interpolation space $\dot{\bf H}^s:=(\dot {\bf L}^2,\dot {\bf H}^2)_{[\frac{s}{2}]}$ for $s\in[1,2]$.

Let the domain $\Omega$ be partitioned into quasi-uniform tetrahedra $K_j$, $j=1,2,\cdots,J$, and denote by $h=\max_j{\rm diam}(K_j)$ the spatial mesh size. 
%
We consider a conforming finite element space ${\bf X}_h\times V_h\subset {\bf H}^1_0\times L^2_0$ with the following approximation properties:
\begin{align}
&\inf_{{\bf v}_h\in {\bf X}_h}\|{\bf v}-{\bf v}_h\|_{L^q}
\le C\|{\bf v}\|_{H^{l+1}}h^{l+1+\frac3q-\frac32} \quad 
&&\forall~{\bf v}\in {\bf H}^1_0\cap {\bf H}^{l+1} , \label{approx1} \\
& &&\mbox{for}\,\, l=0,1\,\,\mbox{and}\,\, 2\le q\le 6,  \nn \\[8pt]
&\inf_{q_h\in V_h}\|q-q_h\|_{L^2}
\le C\|q\|_{H^{1}}h   
&&\forall~q\in H^{1}  , \label{approx2} 
\end{align}
satisfying the inf-sup condition 
\begin{align}\label{inf-sup}
&\|q_h\|_{L^2}\le C\sup_{\begin{subarray}{ll}
{\bf v}_h\in {\bf X}_h\\
{\bf v}_h\neq 0
\end{subarray}}\frac{|(\nabla\cdot{\bf v}_h,q_h)|}{\|{\bf v}_h\|_{H^1}}
\quad
\forall~q_h\in V_h , 
\end{align}
where $C$ is some positive constant independent of the mesh size $h$. In addition, we assume that 
\begin{align}\label{additiona-cond}
&\nabla\cdot \v_h\in V_h\,\,\,\mbox{for}\,\,\, \v_h\in {\bf X}_h ,
\end{align}
which guarantees that the discrete divergence-free functions are divergence-free pointwise, 
a desired property in the numerical solution of the Navier--Stokes equations. Examples of finite element spaces satisfying properties \eqref{approx1}--\eqref{additiona-cond} include the Scott--Vogelius finite element space \cite{2019-Guzman-Scott} and the conforming divergence-free finite element space in \cite{2014-Guzman-Neilan}. 

The inf-sup condition \eqref{inf-sup} guarantees the existence, uniqueness and stability of the finite element solution, but it will not be used explicitly in this paper as we are not going to present error estimates for the pressure. The additional condition \eqref{additiona-cond} is not essential but convenient for error analysis in this paper as it avoids some technical regularity estimates for the pressure in the case of $\dot{\bf H}^2$ initial data.


Let the time interval $[0,T]$ be partitioned uniformly into $0=t_0<t_1<\cdots<t_N=T$, and denote $\tau=t_{m+1}-t_m$. For any sequence of functions $g_0$, $g_1$, ... , $g_N$, we define
$$
D_\tau g^{n+1}=\frac{g^{n+1}-g^n}{\tau } ,
\quad n=0,1,... , N-1 .
$$
For any given $\u_h^n\in {\bf X}_h$, we look for $(\u_h^{n+1},p_h^{n+1})\in {\bf X}_h\times V_h$ as the solution of the following linearized finite element equations
\begin{align}
& \big(D_{\tau}\u^{n+1}_h ,{\bf v}_h\big)
-\big(\u^{n+1}_h,\u^n_h\cdot \nabla {\bf v}_h\big) + \big(\mu \nabla \u^{n+1}_h,\nabla {\bf v}_h\big) 
-\big(p^{n+1}_h,\nabla\cdot{\bf v}_h\big)=0,
\label{FEMEq1}\\[6pt]
&\big(\nabla\cdot \u^{n+1}_h,q_h\big)=0,
\qquad\forall\, {\bf v}_h\in {\bf X}_h\,\,\,\mbox{and}\,\,\, \forall\, q_h\in V_h, 
\label{FEMEq2} 
\end{align}
where $\u_h^0$ is the Stokes--Ritz projection of $\u^0$ onto ${\bf X}_h$ (see section \ref{section:3.2} for the definition of the Stokes--Ritz projection). 

Since the discrete divergence-free subspace coincides with the divergence-free subspace (as a result of \eqref{additiona-cond}), it follows that $\nabla\cdot \u^n_h=0$ and therefore $\big(\u^{n+1}_h,\u^n_h\cdot \nabla \u^{n+1}_h\big) = 0$. As a result, for any given mesh sizes $\tau>0$ and $h>0$, the linearized equations \eqref{FEMEq1}--\eqref{FEMEq2} have a unique finite element solution $\u_h^{n+1}$, $n=0,1,\, ...\, ,N-1$, which satisfy the discrete energy inequality  
\begin{align}\label{EnergEst}
&\max_{0\le n\le N-1}\frac{1}{2}\|\u_h^{n+1}\|_{L^2}^2
+\sum_{n=0}^{N-1}\tau \|\nabla \u_h^{n+1}\|_{L^2}^2
\le \frac{1}{2}\|\u_h^{0}\|_{L^2}^2 \, .
\end{align}

For the solution $\u_h^n$ given by \eqref{FEMEq1}--\eqref{FEMEq2}, we define the piecewise constant numerical solution
\begin{align}\label{DefNSol}
&\u_{h,\tau}(x,t)=\u_h^n(x)\quad\mbox{for}\,\,\, t\in(t_{n-1},t_n] 
\,\,\,\mbox{and}\,\,\, x\in\Omega ,
\end{align}
and present our main result in the following theorem.

\begin{theorem}\label{MainTHM}
{\it
For any $M>0$ there exist positive constants $\tau_M$ and $h_M$ such that 
when 
\begin{align}\label{meshcond0}
\tau<\tau_M
\quad\mbox{and}\quad
h<h_M 
\end{align}
if a numerical solution $\u_{h,\tau}$ defined by \eqref{DefNSol} satisfies
\begin{align}\label{numersolM}
\|\u_{h,\tau}\|_{L^\infty(0,T;L^4)} + \|\u^0\|_{H^2} + 1 \le M,
\end{align}
then the Navier--Stokes equations \eqref{NSPDE1}--\eqref{NSPDE2} possess a unique solution with regularity
\begin{align}\label{RegU}
\u\in L^\infty(0,T;\dot{\bf H}^2) 
\quad\mbox{and}\quad
\partial_t\u\in L^\infty(0,T;\dot{\bf L}^2) .
\end{align}
The constants $\tau_M$ and $h_M$ are decreasing functions of $M$, independent of $\u$, $\u^0$ and $T$, but may depend on $\mu$. 
}
\end{theorem}

\begin{remark}
Theorem \ref{MainTHM} states that, when solving the Navier--Stokes equations, we do not need to assume existence, uniqueness or regularity of the solution. Instead, if we have a initial data $\u^0$ and a numerical solution $\u_{h,\tau}$, one can pick up $M$ satisfying \eqref{numersolM} and refine the mesh according to \eqref{meshcond0}. 
If the conditions \eqref{meshcond0}--\eqref{numersolM} are satisfied by one numerical solution,  then one can say that for the given initial value the Navier--Stokes equations have a unique solution with regularity \eqref{RegU}. 
In this paper we only prove the existence of constants $\tau_M$ and $h_M$. From the expression of $\Phi(M)$ in the proof of Lemma \ref{prop1} and the expression of $\tau_M$ and $h_M$ in \eqref{expr-tauM}--\eqref{expr-hM} one can see that $\tau_M,h_M\ge \exp(-C\mu^{-897}M^{900})$. This estimate for $\tau_M$ and $h_M$ has only theoretical value due to the heavy dependence on $\mu^{-1}$ and $M$, especially for problems with small viscosity (large Reynolds number) and large initial data. 
Nevertheless, there is possibility that more useful estimates may be obtained by different analysis or better error estimates. 
\end{remark}

\begin{remark}
From the proof of the theorem in the next section, one can see that the numerical solution successfully approximates the exact solution in the sense that
\begin{align}\label{ErrEstM}
\|\u_{h,\tau}-\u\|_{L^\infty(0,T;L^2)}^2\le 2(\tau +h^{3/2})  . 
\end{align}
Clearly, the order of this error estimate can be improved. 
The purpose of this paper is to prove the existence of a smooth solution instead of optimal-order error estimates.
\end{remark}

\begin{remark}
The $L^\infty(0,T;L^4)$ norm used in \eqref{numersolM} may be replaced by some other norms stronger than the critical norm of $L^\infty(0,T;L^3)$. 
The analysis can also be extended to the case where the forcing term is not zero, provided that the forcing term is sufficiently smooth in time. In this case, the condition 
\begin{align*}
\|\u_{h,\tau}\|_{L^\infty(0,T;L^4)} + \|\u^0\|_{H^2} + 1 \le M 
\end{align*}
in Theorem \ref{MainTHM} should be replaced by
\begin{align*}
\|\u_{h,\tau}\|_{L^\infty(0,T;L^4)} + \|\u^0\|_{H^2} + \|{\bf f}\|_{L^2(0,T;L^2)} + \|{\bf f}\|_{L^\infty(0,T;L^2)} + \|\partial_t{\bf f}\|_{L^2(0,T;L^2)} + 1 \le M . 
\end{align*}
The three different norms on ${\bf f}$ are all needed as we require the constants $C$ in this paper to be independent of $T$. 
\end{remark}

\begin{remark}
It is possible to extend Theorem \ref{MainTHM} to other nonlinear time-evolution equations for which global existence, uniqueness and regularity of the solution are unknown but local existence, uniqueness and regularity are known for smooth initial data. For such equations, the method could be used to prove global uniqueness and regularity of the solution as well as convergence of the numerical solutions in an a posteriori way. This can be viewed as an improvement of the traditional approach on error estimates of numerical solutions (which is based on global well-posedness assumptions that are not proved yet). 
\end{remark}

\section{Proof of Theorem \ref{MainTHM}}
\setcounter{equation}{0}

It is well known that a solution with the regularity \eqref{RegU} is unique and smooth for $t\in(0,T]$. We shall prove existence of a solution with this regularity up to time $T$. In the rest part of this paper, we denote by $C_{p_1,p_2,...,p_m}$ a generic positive constant which may depend on the parameters $p_1,p_2,...,p_m$ and $\mu$, but is  independent of $n$, $k$, $\tau$, $h$, $T$ and $\mu$. 

Recall that $P:{\bf L}^2\rightarrow \dot{\bf L}^2$ is the $L^2$-orthogonal projection onto the divergence-free subspace. In particular, for any $\v\in {\bf L}^2$ we have $P\v=\v-\nabla q$, with $q\in H^1$ being the weak solution of 
$$
\left\{
\begin{aligned}
&\Delta q=\nabla\cdot\v &&\mbox{in}\,\,\,\Omega,\\
&\nabla q\cdot{\bf n}= \v\cdot{\bf n}  &&\mbox{on}\,\,\,\partial\Omega . 
\end{aligned}
\right. 
$$ 
The $H^2$-regularity estimate of linear Stokes equations in \cite{Dauge-1989} implies that
\begin{align}\label{H2-estimate-Stokes}
\|\v\|_{H^2}
\le
C\|P\Delta\v\|_{L^2} \quad\mbox{for}\,\,\, \v\in \dot {\bf H}^2 . 
\end{align}
The Navier--Stokes equations can be written as 
\begin{align}\label{NS-P}
\partial_t\u+P(\u\cdot\nabla\u)-P\Delta\u=0 . 
\end{align}

\subsection{Local existence and estimates}\label{SEction3}
In this subsection, we prove the following lemma, which is used in the next subsection to prove Theorem \ref{MainTHM}.
\begin{lemma}\label{prop1}
There exists a decreasing function $\varphi:\R_+\rightarrow\R_+$ and an increasing function $\Phi:\R_+\rightarrow\R_+$, with $\Phi(s)\ge s$,  such that if $\u^0\in \dot{\bf H}^2$ and the solution of {\rm(\ref{NSPDE1})--(\ref{NSPDE2})} has the regularity \eqref{RegU} up to time $T$, then the solution can be extended to time $T+\varphi(\|\u\|_{L^\infty(0,T;L^4)} +\|\u^0\|_{H^2})$ with the same regularity, satisfying the following quantitative estimate: 
\begin{align*}
&\|\partial_{tt}\u\|_{L^2(0,T+\varphi(\|\u\|_{L^\infty(0,T;L^4)}
+\|\u^0\|_{H^2}) ;\dot H^{-1})} 
+\|\partial_{t}\u\|_{L^\infty(0,T+\varphi(\|\u\|_{L^\infty(0,T;L^4)}
+\|\u^0\|_{H^2}) ;L^2)} \nn\\
&\quad
+\|\partial_{t}\u\|_{L^2(0,T+\varphi(\|\u\|_{L^\infty(0,T;L^4)}
+\|\u^0\|_{H^2}) ;H^1)}
+ \|\u\|_{L^2(0,T+\varphi(\|\u\|_{L^\infty(0,T;L^4)}
+\|\u^0\|_{H^2});H^2)} \nn\\
&\quad
+ \|\u\|_{L^\infty(0,T+\varphi(\|\u\|_{L^\infty(0,T;L^4)}
+\|\u^0\|_{H^2});H^2)}
+ \|p\|_{L^\infty(0,T+\varphi(\|\u\|_{L^\infty(0,T;L^4)}
+\|\u^0\|_{H^2});H^1)}  \nn\\
&\le \Phi(\|\u\|_{L^\infty(0,T;L^4)}
+\|\u^0\|_{H^2})   ,
\end{align*}
where the functions $\varphi$ and $\Phi$ do not depend on $\u$ or $T$.

\end{lemma}


\begin{remark}
The property $\Phi(s)\ge s$ is not necessarily needed. It is used only to simplify the notation in the proof of Theorem \ref{MainTHM}. 
\end{remark}

\begin{proof}

In order to prove Lemma \ref{prop1}, we introduce some lemmas below.

\begin{lemma}\label{Lem0}
{\it In a three-dimensional convex domain $\Omega$, there exists a positive constant $C_0$ such that
\begin{align*}
&\|\v\|_{L^p} \le C_0\|\v\|_{H^1} 
&& \forall\, 1\le p\le 6
&&\forall\,\v\in {\bf H}^1  \\
&\|\v\|_{L^p} \le C_0\|\v\|_{L^2}^{1-\theta}\|\nabla\v\|_{L^2}^{\theta}
&&\mbox{with}\,\,\,p\in[2,6]\,\,\,\mbox{and}\,\,\,\,\frac{1-\theta}{2}+\frac{\theta}{6} = \frac{1}{p}  &&\forall\,\v\in {\bf H}^1_0  \\
&\|\v\|_{L^4} \le C_0\|\v\|_{L^2}^{\frac14}\|\nabla\v\|_{L^2}^{\frac34} &&
&&\forall\,\v\in {\bf H}^1_0  \\
&\|\nabla\v\|_{L^4} \le C_0\|\nabla\v\|_{L^2}^{\frac14}\|\v\|_{H^2}^{\frac34} 
&&
&&\forall\,\v\in {\bf H}^2 \\
&
\|\v\|_{L^\infty} \le C_0\|\v\|_{H^1}^{\frac12}\|\v\|_{H^2}^{\frac12}   
&&
&& \forall\,\v\in {\bf H}^2 
\end{align*}
}
\end{lemma}

\begin{lemma}\label{Lem2}
{\it There exists an increasing function $\beta:\R_+\rightarrow \R_+$ such that if $\u^0\in \dot{\bf H}^2$ and if the Navier--Stokes equations \eqref{NSPDE1}--\eqref{NSPDE2} has a weak solution $\u\in L^\infty(0,T;\dot{\bf H}^1)$, then the solution has regularity \eqref{RegU} and satisfies the following quantitative estimate:  
\begin{align*}
&\|\partial_{tt}\u\|_{L^2(0,T;\dot H^{-1})} 
+\|\partial_{t}\u\|_{L^\infty(0,T;L^2)} +\|\partial_{t}\u\|_{L^2(0,T;H^1)} \nn\\
&\quad
+\|\u\|_{L^2(0,T;H^2)} + \|\u\|_{L^\infty(0,T;H^2)} + \|p\|_{L^\infty(0,T;H^1)} 
\nn\\ 
&\le \beta(\|\u\|_{L^\infty(0,T;L^4)}+\|\u^0\|_{H^2}) ,
\end{align*}
where the function $\beta$ does not depend on $\u$ or $T$. 
}
\end{lemma}

\begin{lemma}\label{Lem1}
{\it There exists a decreasing function $\alpha:\R_+\rightarrow \R_+$ such that if $\u^0\in \dot{\bf H}^2$, then the Navier--Stokes equations  {\rm(\ref{NSPDE1})--(\ref{NSPDE2}) have a unique strong solution with regularity \eqref{RegU} up to time $T=\alpha(\|\u^0\|_{H^2})$, satisfying the following quantitative estimate:} 
$$
\|\u\|_{L^\infty(0,\alpha(\|\u^0\|_{H^2});H^1)}
\le  \|\u^0\|_{H^1} +1 ,
$$
where the function $\alpha$ does not depend on $\u$ or $T$. 
}
\end{lemma}

\begin{remark}
The first four inequalities in Lemma \ref{Lem0} are consequences of the Sobolev embedding inequality \cite[p. 102, Theorem 4.31]{Adams} and the Sobolev interpolation inequality \cite[p. 139, Theorem 5.8]{Adams}.
The last inequality in Lemma \ref{Lem0} is Agmon's inequality. 
The proof of Agmon's inequality in a bounded smooth domain can be found in \cite[Lemma 4.10]{Constantin-Foias}. In a general bounded Lipschitz domain, there exists a linear extension operator $E: {\bf L}^1(\Omega)\rightarrow {\bf L}^1(\R^3)$, e.g., Stein's extension operator in \cite[p.\ 181, Theorem 5]{Stein1970}, which satisfies that  
\begin{enumerate}
\item
$E{\bf v}={\bf v}$ a.e. in $\Omega$.

\item
$E{\bf v}$ is supported in a ball $B$ containing $\overline\Omega$ (this can be ensured by multiplying the extended functions with a common smooth cut-off function that is equal to $1$ on $\overline\Omega$ and equal to $0$ outside $B$).

\item
For all $k\ge 0$ and $1\le p\le \infty$, if ${\bf v}\in W^{k,p}(\Omega)$ then 
$\|E{\bf v}\|_{W^{k,p}(\R^3)}\le C_{k,p}\|{\bf v}\|_{W^{k,p}(\Omega)}$.  
\end{enumerate}
As a result, 
$$
\|\v\|_{L^\infty(\Omega)}
\le
\|E\v\|_{L^\infty(B)} 
\le C\|E\v\|_{H^1(B)}^{\frac12}\|E\v\|_{H^2(B)}^{\frac12} 
\le
C\|\v\|_{H^1(\Omega)}^{\frac12}\|\v\|_{H^2(\Omega)}^{\frac12} . 
$$
This verifies that Agmon's inequality holds on general bounded Lipschitz domains. 

Based on the proof of Lemmas \ref{Lem2}--\ref{Lem1} below, one can choose 
$$
\alpha(s)=\frac{1}{[C_1^*\mu^{-\frac{1}{2}}+C_1^*\mu^{-\frac{15}{2}}[C_0+(C_0+1)s]^{9}]^2}
\quad\mbox{and}\quad
\beta(s)=C_1\mu^{-2}+C_1\mu^{-14}s^{15} ,
$$ 
where $C_1$ and $C_1^*$ are some positive constants independent of $\u$ and $T$.
\end{remark}

Lemma \ref{prop1} assumes that the solution is a strong solution with regularity \eqref{RegU}, i.e., $\u\in L^\infty(0,T;\dot{\bf H}^2) \cap W^{1,\infty}(0,T;\dot{\bf L}^2) $. 
For $s\in(\frac32,2)$ the following interpolation inequality (\cite[p. 59, Theorem (f)]{Triebel}) is known: 
\begin{align*}
\|\u(t)-\u(t')\|_{\dot H^s}
&\le
C\|\u(t)-\u(t')\|_{\dot L^2}^{1-\frac{s}{2}}\|\u(t)-\u(t')\|_{\dot H^2}^{\frac{s}{2}} \\
&\le
C|t-t'|^{1-\frac{s}{2}}
\|\u\|_{W^{1,\infty}(0,T;L^2)}^{1-\frac{s}{2}} \|\u\|_{L^\infty(0,T;H^2)}^{\frac{s}{2}} ,
\end{align*}
which implies that $L^\infty(0,T;\dot{\bf H}^2) \cap W^{1,\infty}(0,T;\dot{\bf L}^2) \hookrightarrow C([0,T];\dot{\bf H}^s)$ for $s\in(\frac32,2)$. Hence, $\u(\cdot,T)\in \dot{\bf H}^s$. 
Since $\dot{\bf H}^s\hookrightarrow {\bf L}^\infty $ for $s\in(\frac32,2)$, based on Taylor's result \cite[Proposition 1.1]{Taylor}, the solution can be furthermore extended to $C([0,T+\epsilon];\dot{\bf H}^s)\hookrightarrow L^\infty(0,T+\epsilon;\dot{\bf L}^6)$. This satisfies Serrin's regularity condition (cf. \cite[section 3, (15)--(16)]{Serrin}) in the time interval $(0,T+\epsilon)$. In this case, the solution is qualitatively $C^\infty(\overline\Omega_K\times[\frac{T}{2},T])$ on any subset $\Omega_K$ such that $\overline\Omega_K\subset\Omega$; cf. \cite[section 3, (15)--(16)]{Serrin}. In particular, 
$$
\|\u(\cdot,T)\|_{H^2(\Omega_K)}
\le 
C\|\u\|_{C([\frac{T}{2},T]; H^2(\Omega_K))} 
\le 
C\|\u\|_{L^\infty(0,T; H^2(\Omega))} ,
$$
where the right-hand side is bounded due to the assumption in Lemma \ref{prop1} that the solution has regularity \eqref{RegU}. 
Since the right-hand side above is independent of the subset $\Omega_K$, it follows that $\u(\cdot,T)\in \dot{\bf H}^2$. Then Lemma \ref{Lem1} implies that the solution can be furthermore extended to time $T+\alpha(\|\u(\cdot,T)\|_{H^2})$, satisfying the qualitative regularity 
$$
\u\in L^\infty(0,T+\alpha(\|\u(\cdot,T)\|_{H^2});\dot {\bf H}^{2})\cap 
W^{1,\infty}(0,T+\alpha(\|\u(\cdot,T)\|_{H^2});\dot {\bf L}^{2}) 
$$ 
and the quantitative estimate: 
\begin{align*} 
& \|\u\|_{L^\infty( 0,T+\alpha(\|\u(\cdot,T)\|_{H^2}); H^1)} 
\le \|\u\|_{L^\infty(0,T;H^1)}+1 . 
\end{align*} 
Since $\alpha(\cdot)$ is a decreasing function and $\|\u(\cdot,T)\|_{H^2}\le \|\u\|_{L^\infty(0,T;H^2)}$, it follows that
\linebreak $T+\alpha(\|\u\|_{L^\infty(0,T;H^2)})\le T+\alpha(\|\u(\cdot,T)\|_{H^2})$. Hence, $\u\in L^\infty(0,T+\alpha(\|\u\|_{L^\infty(0,T;H^2)});\dot{\bf H}^2)$ and 
\begin{align*}
& \|\u\|_{L^\infty( 0,T+\alpha(\|\u\|_{L^\infty(0,T;H^2)}); H^1)} 
\le \|\u\|_{L^\infty(0,T;H^1)}+1 .
\end{align*}

By considering the solution in the time interval $(0,T)$ in Lemma \ref{Lem2}, we obtain
\begin{align*}
\|\u\|_{L^\infty(0,T;H^2)}
\le
\beta(\|\u\|_{L^\infty(0,T;L^4)}+\|\u^0\|_{H^2}) .
\end{align*}
By considering the solution in the time interval $(0,T+\alpha(\|\u\|_{L^\infty(0,T;H^2)}))$ in Lemma \ref{Lem2}, we obtain 
\begin{align*}
&\|\partial_{tt}\u\|_{L^2(0,T+\alpha(\|\u\|_{L^\infty(0,T;H^2)});\dot H^{-1})} 
+\|\partial_{t}\u\|_{L^\infty(0,T+\alpha(\|\u\|_{L^\infty(0,T;H^2)});L^2)} \\
&\quad 
+\|\partial_{t}\u\|_{L^2(0,T+\alpha(\|\u\|_{L^\infty(0,T;H^2)});H^1)}
+ \|\u\|_{L^2(0,T+\alpha(\|\u\|_{L^\infty(0,T;H^2)});H^2)} \\
&\quad
+ \|\u\|_{L^\infty(0,T+\alpha(\|\u\|_{L^\infty(0,T;H^2)});H^2)}+ \|p\|_{L^\infty(0,T+\alpha(\|\u\|_{L^\infty(0,T;H^2)});H^1)} \\
&\le \beta(\|\u\|_{L^\infty(0,T+\alpha(\|\u\|_{L^\infty(0,T;H^2)});L^4)}+\|\u^0\|_{H^2}) \\ 
&\le \beta(C_0\|\u\|_{L^\infty(0,T+\alpha(\|\u\|_{L^\infty(0,T;H^2)});H^1)}+\|\u^0\|_{H^2})  ,
\end{align*}
where $C_0$ is the constant in Lemma \ref{Lem0}. The last three estimates imply Lemma \ref{prop1} with
$$
\varphi(s)=\alpha(\beta(s))\qquad\mbox{and}\qquad
\Phi(s)=\beta(C_0\beta(s)+C_0+s) . 
$$

It remains to prove Lemma \ref{Lem2} and Lemma \ref{Lem1}.
\end{proof}

\begin{proof}[Proof of Lemma \ref{Lem2}]
Under the conditions of Lemma \ref{Lem2} we have $\u\in L^\infty(0,T;\dot {\bf H}^1)\hookrightarrow L^\infty(0,T;{\bf L}^6)$, which satisfies Serrin's regularity condition (cf. \cite[section 3, (15)--(16)]{Serrin}). In this case, the solution to the Navier--Stokes equations is qualitatively smooth in the domain $\Omega\times(0,T)$. In the following, we present quantitative estimates for the solution with all the positive constants independent of $\u$ and $T$. 

First, integrating \eqref{NSPDE1} against $\u$ yields the basic energy estimate 
\begin{align}\label{basic-energy}
\frac12 \|\u\|_{L^\infty(0,T;L^2)}^2
+\mu\|\nabla \u\|_{L^2(0,T;L^2)}^2 
&\le \frac12 \|\u^0\|_{L^2}^2 .
\end{align}

Second, integrating \eqref{NSPDE1} against $\partial_t\u$ yields 
\begin{align*}
\|\partial_t\u\|_{L^2}^2
+\frac{\d}{\d t}\bigg(\frac{\mu}{2}\|\nabla \u\|_{L^2}^2\bigg)
&\leq\frac{1}{2}\|\u\cdot\nabla\u\|_{L^2}^2
+\frac{1}{2}\|\partial_t\u\|_{L^2}^2  ,
\end{align*}
which can be reduced to 
\begin{align*}
\|\partial_t\u\|_{L^2}^2
+\frac{\d}{\d t}\bigg( \mu \|\nabla \u\|_{L^2}^2\bigg)
&\le \|\u\cdot\nabla\u\|_{L^2}^2 .
\end{align*}
From \eqref{NS-P} we furthermore derive that 
\begin{align*}
\mu^2 \|P\Delta\u\|_{L^2}^2
&\le 2\|\partial_t\u\|_{L^2}^2
+ 2\|P(\u\cdot\nabla\u)\|_{L^2}^2 \\
&\le 4\|\u\cdot\nabla\u\|_{L^2}^2   
-\frac{\d}{\d t}\bigg( 2\mu \|\nabla \u\|_{L^2}^2\bigg).
\end{align*}
The sum of the last two estimates gives 
\begin{align*}
&\|\partial_t\u\|_{L^2}^2
+ \mu^2 \|P\Delta\u\|_{L^2}^2
+\frac{\d}{\d t}\bigg( 3\mu \|\nabla \u\|_{L^2}^2\bigg) \\
&\le 5\|\u\|_{L^4}^2\|\nabla\u\|_{L^4}^2 \\
&\le C\|\u\|_{L^4}^{2}\|\nabla\u\|_{L^2}^{\frac12}\|\u\|_{H^2}^{\frac32} \\
&\le C\mu^{-\frac32}\|\u\|_{L^4}^{2}\|\nabla\u\|_{L^2}^{\frac12} (\mu^{\frac32}\|P\Delta\u\|_{L^2}^{\frac32}) \\
&\le C\mu^{-6} \|\u\|_{L^4}^{8}\|\nabla\u\|_{L^2}^{2}
+ \frac{\mu^2}{2}\|P\Delta\u\|_{L^2}^2 ,
\end{align*} 
where we have used Lemma \ref{Lem0} in the third to last inequality and \eqref{H2-estimate-Stokes} in the second to last inequality. Replacing $\|P\Delta\u\|_{L^2}^2$ by $C\|\u\|_{H^2}^2$ on the left-hand side and integrating the result in time, we obtain
\begin{align}\label{L2H2u}
&\|\partial_t\u\|_{L^2(0,T;L^2)}^2
+C\mu^2 \|\u\|_{L^2(0,T;H^2)}^2
+  3\mu \|\nabla \u\|_{L^\infty(0,T;L^2)}^2 \nn \\
&\le 3\mu \|\nabla \u^0\|_{L^2}^2
+C\mu^{-6}\|\u\|_{L^\infty(0,T;L^4)}^{8}\|\nabla\u\|_{L^2(0,T;L^2)}^{2} \nn \\
&\le 3\mu \|\nabla \u^0\|_{L^2}^2
+C\mu^{-7}\|\u\|_{L^\infty(0,T;L^4)}^{8}\|\u^0\|_{L^2}^{2} \nn\\
&\le C\mu+C\mu^{-7}(\|\u\|_{L^\infty(0,T;L^4)}+\|\u^0\|_{H^2})^{10} ,
\end{align}
where we have used \eqref{basic-energy} in obtaining the second to last inequality and again, Young's inequality in deriving the last inequality. 

Then, differentiating \eqref{NSPDE1} with respect to $t$, we have 
\begin{align}
&\partial_{tt}\u-\mu\Delta \partial_t\u+\nabla \partial_tp
=-\nabla\cdot(\partial_t\u\otimes\u)-\nabla\cdot(\u\otimes\partial_t\u) .
\label{Equtt}
\end{align}
Then, integrating the equation above against $\mu\partial_t\u$, we obtain
\begin{align*}
\frac{\d}{\d t}\left(\frac{\mu}{2}\|\partial_{t}\u\|_{L^2}^2\right)
+ \mu^2 \|\nabla\partial_{t}\u\|_{L^2}^2
&\le 2\|\partial_t\u\otimes\u\|_{L^2}^2 + \frac{\mu^2}{2} \|\nabla\partial_{t}\u\|_{L^2}^2 
\end{align*}
and therefore
\begin{align*}
\frac{\d}{\d t}\left(\frac{\mu}{2}\|\partial_{t}\u\|_{L^2}^2\right)
+\frac{\mu^2}{2}\|\nabla\partial_{t}\u\|_{L^2}^2
&\le 2\|\partial_t\u\otimes\u\|_{L^2}^2 \\
&\le  C\|\partial_t\u\|_{ L^4 }^2
\|\u\|_{L^4}^2 \\
&\le C\|\partial_t\u\|_{L^2 }^{\frac12}
\|\partial_t\u\|_{L^6}^{\frac32}
\|\u\|_{L^4}^2 \\
&\le C\|\partial_t\u\|_{L^2 }^{\frac12}
\|\nabla\partial_t\u\|_{L^2 }^{\frac32}
\|\u\|_{L^4}^2\\
&\le C\mu^{-6}\|\u\|_{L^4}^8\|\partial_t\u\|_{L^2 }^2
+ \frac{\mu^2}{4} \|\nabla\partial_t\u\|_{L^2 }^{2} ,
\end{align*}
where we have used the Sobolev embedding ${\bf H}^1\hookrightarrow {\bf L}^6$ in the second to last inequality. After absorbing the last term of the inequality above by its left-hand side, integrating the result in time yields 
\begin{align}\label{uttL2H-1}
&2 \mu \|\partial_{t}\u\|_{L^\infty(0,T;L^2)}^2 
+ \mu^2 \|\nabla\partial_{t}\u\|_{L^2(0,T;L^2)}^2 \nn\\
&\le 2\mu \|\partial_{t}\u^0\|_{L^2}^2+C\mu^{-6} \|\u\|_{L^\infty(0,T;L^4)}^8
\|\partial_t\u\|_{L^2(0,T;L^2) }^2 \nn\\
&\le 
2\mu \|\u^0\cdot\nabla\u^0-\mu\Delta \u^0\|_{L^2}^2
\quad\mbox{(here we used $\partial_t\u^0=P[\u^0\cdot\nabla\u^0-\mu\Delta \u^0]$)}  \nn\\
&\quad\, 
+C \mu^{-6} \|\u\|_{L^\infty(0,T;L^4)}^8[C\mu+C\mu^{-6}( 
\|\u\|_{L^\infty(0,T;L^4)} +\|\u^0\|_{H^2} )^{10} ] \nn\\
&\le C\mu +C\mu^{-13}(\|\u\|_{L^\infty(0,T;L^4)}+\|\u^0\|_{H^2})^{18}  , 
\end{align}
where we have used \eqref{L2H2u} in estimating $ \|\partial_{t}\u\|_{L^2(0,T;L^2)}^2$.  
Substituting \eqref{uttL2H-1} into \eqref{Equtt} and using the duality argument (testing the equation by a function in $\dot{\bf H}^1$), we can obtain that
\begin{align}\label{uttL2H-2} 
\|\partial_{tt}\u\|_{L^2(0,T;\dot H^{-1})} 
&\le 
C\mu\|\nabla\partial_{t}\u\|_{L^2(0,T;L^2)} 
+C\|\partial_t\u\otimes \u\|_{L^2(0,T;L^2)}  \nn\\
&\le C\mu^{\frac{1}{2}}+C\mu^{-\frac{13}{2}}(\|\u\|_{L^\infty(0,T;L^4)}+\|\u^0\|_{H^2})^{9}   .
\end{align}

Next, from \eqref{NSPDE1} and the basic $H^2$ estimate of Stokes equations we know that 
\begin{align*}
\|\u\|_{L^\infty(0,T;H^2)} 
&\le C\|P\Delta \u\|_{L^\infty(0,T;L^2)} \\
&\le C\mu^{-1} \|\partial_t\u\|_{L^\infty(0,T;L^2)}
+C\mu^{-1} \|\u\cdot\nabla\u\|_{L^\infty(0,T;L^2)}\\
&\le C\mu^{-1} \|\partial_t\u\|_{L^\infty(0,T;L^2)}
+ C\mu^{-1} \|\u\|_{L^\infty(0,T;L^\infty)}\|\nabla\u\|_{L^\infty(0,T;L^2)} \\
&\le C\mu^{-1} \|\partial_t\u\|_{L^\infty(0,T;L^2)}+ C\mu^{-1} \|\u\|_{L^\infty(0,T;H^1)}^{\frac12}\|\u\|_{L^\infty(0,T;H^2)}^{\frac12}
\|\nabla\u\|_{L^\infty(0,T;L^2)} \\
&\le C\mu^{-1} \|\partial_t\u\|_{L^\infty(0,T;L^2)}
+ C\mu^{-2} \|\nabla\u\|_{L^\infty(0,T;L^2)}^3
+ \frac{1}{2} \|\u\|_{L^\infty(0,T;H^2)} ,
\end{align*}
where we have used the last inequality of Lemma \ref{Lem0} (Agmon's inequality) in deriving the second to last inequality.  
The inequality above implies 
\begin{align}\label{LinfH2u}
\|\u\|_{L^\infty(0,T;H^2)}
&\le C\mu^{-1} \|\partial_t\u\|_{L^\infty(0,T;L^2)}
+ C\mu^{-2} \|\nabla\u\|_{L^\infty(0,T;L^2)}^3 \nn \\
&\le 
C\mu^{-1}+C\mu^{-8}(\|\u\|_{L^\infty(0,T;L^4)}+\|\u^0\|_{H^2})^{9} \nn \\
&\quad\, + C\mu^{-2}+C\mu^{-14}(\|\u\|_{L^\infty(0,T;L^4)}+\|\u^0\|_{H^2})^{15} \nn \\
&\le 
C\mu^{-2}+C\mu^{-14}(\|\u\|_{L^\infty(0,T;L^4)}+\|\u^0\|_{H^2})^{15} , 
\end{align}
where we have used \eqref{L2H2u} and \eqref{uttL2H-1} in estimating $\|\nabla\u\|_{L^\infty(0,T;L^2)}$ and $\|\partial_t\u\|_{L^\infty(0,T;L^2)}$, respectively. 
With the above estimates of $\|\partial_t\u\|_{L^\infty(0,T;L^2)}$,
$ \|\u\cdot\nabla\u\|_{L^\infty(0,T;L^2)}$
and $\|\u\|_{L^\infty(0,T;H^2)}$, from \eqref{NSPDE1} we also derive that
\begin{align}\label{LinfH1p}
\|p\|_{L^\infty(0,T;H^1)}
&\le C\mu^{-1}+C\mu^{-13}(\|\u\|_{L^\infty(0,T;L^4)}+\|\u^0\|_{H^2})^{15} .
\end{align}

Finally, the inequalities \eqref{L2H2u} and \eqref{uttL2H-1}--\eqref{LinfH1p}
imply Lemma \ref{Lem2}
with 
\begin{align} 
\beta(s)=C_1\mu^{-2}+C_1\mu^{-14}s^{15} 
\end{align}
where $C_1$ is some positive constant independent of $\u$ and $T$. Without loss of generality, we can choose the constant $C_1$ to be bigger than $1$ so that $\beta(s)\ge s$ for $s\ge 0$.
\end{proof}

\noindent{\it Proof of Lemma \ref{Lem1}.}$\quad$
Let $s\in(\frac32,2)$ be a fixed number. Then $\dot{\bf H}^2\hookrightarrow \dot{\bf H}^s\hookrightarrow {\bf L}^\infty$. In this case, Taylor's local existence result \cite[Proposition 1.1]{Taylor} says that for a given initial value $\u^0\in \dot{\bf H}^2\hookrightarrow \dot{\bf H}^s$, there exists $T_0>0$ such that the Navier--Stokes equations have a unique weak solution in $C([0,T_0];\dot{\bf H}^s)$. We denote by $T_*$ the supremum of such $T_0$, namely,  
$$
\u\in C([0,T_0];\dot{\bf H}^s)\quad\forall\, T_0\in(0,T_*)  
\quad\mbox{and}\quad
\u\notin C([0,T_*];\dot{\bf H}^s) .
$$ 

If we denote by $t_*$ the supremum of time $t>0$ such that the Navier--Stokes equations with initial value $\u^0$ have a weak solution $\u\in L^\infty(0,t;\dot{\bf H}^1)$ satisfying 
$$
\|\u\|_{L^\infty(0,t;H^1)} 
\le \|\u^0\|_{H^1}+1 ,
$$ 
then $\|\u\|_{L^\infty(0,t_*;H^1)} \le \|\u^0\|_{H^1}+1$ and $\u\in L^\infty(0,t_*;\dot{\bf H}^1)$ is a weak solution satisfying the conditions of Lemma \ref{Lem2}, which implies that the regularity of the solution can be furthermore picked up to 
\begin{align} \label{Linfty-H2-t-star}
\u 
&\in L^\infty(0,t_*;\dot{\bf H}^2)\cap W^{1,\infty}(0,t_*;\dot{\bf L}^2) \\
& \hookrightarrow C([0,t_*];\dot{\bf H}^s)\quad\mbox{for}\,\,\,s\in(\mbox{$\frac32$},2). \nonumber 
\end{align} 
This implies $t_*<T_*$. 

The regularity $\u\in C([0,t_*];\dot{\bf H}^s)\hookrightarrow C([0,t_*];\dot{\bf H}^1)$ and the definition of $t_*$ imply 
$$\|\u^0\|_{H^1}+1=\|\u\|_{L^\infty(0,t_*;H^1)} .$$
By using the Newton--Leibnitz formula and the estimate of $\|\partial_t\u\|_{L^2(0,t_*;H^1)}$ in \eqref{uttL2H-1}, we derive that
\begin{align*} 
&\|\u^0\|_{H^1}+1=\|\u\|_{L^\infty(0,t_*;H^1)} \\ 
&\le \|\u^0\|_{H^1}+\|\partial_t\u\|_{L^2(0,t_*;H^1)}t_*^{1/2} \\ 
&\le \|\u^0\|_{H^1}
+[C_1^*\mu^{-\frac12}+C_1^*\mu^{-\frac{15}{2}}(\|\u\|_{L^\infty(0,t_*;L^4)}
+\|\u^0\|_{H^2})^{9}]t_*^{1/2}\\
&\le \|\u^0\|_{H^1}
+[C_1^*\mu^{-\frac{1}{2}}+C_1^*\mu^{-\frac{15}{2}}(C_0\|\u\|_{L^\infty(0,t_*;H^1)}
+\|\u^0\|_{H^2})^{9}]t_*^{1/2}  ,
\end{align*}
which implies 
\begin{align*} 
t_* 
&\geq \frac{1}{[C_1^*\mu^{-\frac{1}{2}}+C_1^*\mu^{-\frac{15}{2}}(C_0+C_0\|\u^0\|_{H^1}+\|\u^0\|_{H^2})^{9}]^2} \\ 
&\geq \frac{1}{[C_1^*\mu^{-\frac{1}{2}}+C_1^*\mu^{-\frac{15}{2}}[C_0+(C_0+1)\|\u^0\|_{H^2}]^{9}]^2} .
\end{align*}

In view of \eqref{Linfty-H2-t-star} and the above lower bound of $t_*$, Lemma \ref{Lem1} holds with 
$$\alpha(s)=\frac{1}{[C_1^*\mu^{-\frac{1}{2}}+C_1^*\mu^{-\frac{15}{2}}[C_0+(C_0+1)s]^{9}]^2} .$$
\qed\medskip

\subsection{Global existence and estimates based on a numerical solution}
\label{section:3.2}

We introduce the Stokes--Ritz projection operator
$(R_h,P_h): {\bf H}^1_0\times L^2  
\rightarrow {\bf X}_h\times V_h$ by
\begin{align}\label{RitzP}
\begin{aligned}
&\big(\nabla({\bf w}-R_h({\bf w},p)),\nabla{\bf v}_h\big)
-\big(p-P_h({\bf w},p) , \nabla \cdot {\bf v}_h\big) =0,
\quad\forall\,{\bf v}_h\in {\bf X}_h,\\
&\big(\nabla\cdot R_h({\bf w},p),q_h\big)=0,
\quad\forall\, q_h\in   V_h ,
\end{aligned}
\end{align} 
and impose the condition $\int_\Omega (p-P_h({\bf w},p))\d x=0$ for uniqueness. 
This Stokes--Ritz projection has the following approximation property: 
\begin{align}\label{RitzPEst0}
&h^{\frac32-\frac3q}\|{\bf w}-R_h({\bf w},p)\|_{L^q}
+h\|{\bf w}-R_h({\bf w},p)\|_{H^1}
+h\|p-P_h({\bf w},p)\|_{L^2} \nn\\
&\le Ch^{l+1}
(\| {\bf w}\|_{H^{l+1}}+\| p\|_{H^{l}}) ,\quad l=0,1 ,\quad
2\le q\le 6 ,\,\,\,
\forall\,
({\bf w},p)\in \dot{\bf H}^2\times H^1,
\end{align}
see \cite{Verfurth} for the proof of the case $q=2$; the case $2< q\le 6$ can be obtained by using the inverse inequality and the Bramble--Hilbert lemma. 

Let $\mathring {\bf X}_h$ be the divergence-free subspace of ${\bf X}_h$, which coincides with the  divergence free subspace of ${\bf X}_h$ undercondition \eqref{additiona-cond}. 
Then the Ritz projection defined in \eqref{RitzP} satisfies that $R_h({\bf w},p)\in \mathring {\bf X}_h$ and $P_h({\bf w},p) \in V_h$, and 
\begin{align}\label{error-Ritz-2}
\begin{aligned}
&\big(\nabla({\bf w}-R_h({\bf w},p)),\nabla{\bf v}_h\big) =0,
&&\forall\,{\bf v}_h\in \mathring {\bf X}_h,\\
&\big(p-P_h({\bf w},p) , \nabla \cdot {\bf v}_h\big) = \big(\nabla({\bf w}-R_h({\bf w},p)),\nabla{\bf v}_h\big),
&&\forall\, {\bf v}_h\in {\bf X}_h ,
\end{aligned}
\end{align} 
Therefore, the operators $R_h$ and $P_h$ are decoupled. In this case, $R_h({\bf w},p)$ is independent of $p$ and therefore error estimate \eqref{RitzPEst0} can be changed to 
\begin{align}\label{RitzPEst2}
&h^{\frac32-\frac3q} \|{\bf w}-R_h({\bf w},p)\|_{L^q}
+h\|{\bf w}-R_h({\bf w},p)\|_{H^1} \nn\\
&\le Ch^{l+1} \| {\bf w}\|_{H^{l+1}} ,\quad l=0,1 ,\quad
2\le q\le 6 ,\,\,\,
\forall\,
({\bf w},p)\in \dot{\bf H}^2\times H^1.
\end{align}

For the simplicity of notation, we denote $R_h{\bf u}=R_h({\bf u},p)$ (as it is independent of $p$) and 
\begin{align}\label{Ritz-simple}
(\u_{*,h},p_{*,h})=(R_h{\bf u},P_h({\bf u},p))
\end{align} 
in the rest of this paper. The error bound in \eqref{error-Ritz-2} implies that  
\begin{align}\label{RitzPEst}
&h^{\frac32-\frac3q}\|{\bf u}-\u_{*,h}\|_{L^q}
+h\|{\bf u}-\u_{*,h}\|_{H^1}
\le Ch^{l+1}
\| {\bf u}\|_{H^{l+1}}  ,\quad l=0,1 ,\,\,\,
2\le q\le 6 .
\end{align}
This approximation property, together with the inverse inequality 
\begin{align}\label{InvInEq}
\|{\bf v}_h\|_{W^{1,q_2}}\le Ch^{\frac{3}{q_2} - \frac{3}{q_1} +l-1}\|{\bf v}_h\|_{W^{l,q_1}} 
,\,\,\,\, \forall\,{\bf v}_h\in {\bf X}_h,\,\,\,
1\le q_1\le q_2\leq\infty,\,\,\, l=0,1 ,
\end{align} 
will be used in the following analysis. We also need
the following version of discrete Gronwall's inequality (cf. \cite{HR90}).
\begin{lemma}
\label{gronwall}
{\it 
Let $\tau$, $B$ and $a_{m}$, $b_{m}$, $c_{m}$, $\gamma_{m}$, 
for integers $k \geq 0$, be nonnegative numbers such that
\[
a_{n+1} + \tau \sum_{m=0}^{n} b_{m+1} 
\le \tau \sum_{m=0}^{n} \gamma_{m} a_{m} + 
\tau \sum_{m=0}^{n} c_{m+1} + B \, , \quad \mathrm{for } \quad n \geq 0 \, .
\]
Then
\[
a_{n+1} + \tau \sum_{m=0}^{n} b_{m+1} 
\le  \exp\bigg(\sum_{m=0}^{n} 
\tau  \gamma_{m} \bigg) \bigg(\tau \sum_{m=0}^{n} c_{m+1} + B\bigg) \, , 
\quad \mathrm{for } \quad n \geq 0 \, .
\]
}
\end{lemma}

To prove existence of a strong solution up to time $T$, we use mathematical induction on $k$ by assuming that 
\begin{align}\label{mathind}
\begin{aligned}
&\mbox{problem \eqref{NSPDE1}--\eqref{NSPDE2} with initial value $\u^0$ has a unique strong solution } \\
&\mbox{$\u\in L^\infty(0,t_k;\dot{\bf H}^2)\cap W^{1,\infty}(0,t_k;\dot{\bf L}^2)$ satisfying $\|\u_{h,\tau}-\u\|_{L^\infty(0,t_k;L^4)}\le 1$}.
\end{aligned}
\end{align}
In the following, we prove that if \eqref{mathind} holds for some nonnegative integer $0\le k\le N-1$, then it also holds when $t_k$ is replaced by $t_{k+1}$. 

To simplify the notation, we denote 
$$
M=\| \u_{h,\tau}\|_{L^\infty(0,T;L^4)}+1+\|\u^0\|_{H^2} .
$$ 
Since $M\ge 1$ and the function $\Phi$ in Lemma \ref{prop1} satisfies $\Phi(M)\ge M$, it follows that $$\Phi(M)\ge 1  .$$ 

Since $\u_h^0$ is the Stokes--Ritz projection of $\u^0$, we have $\|\u_h^0-\u^0\|_{L^4}\le C_2\|\u^0\|_{H^2}h^{5/4}$ for some positive constant $C_2$. Thus our assumption holds for $k=0$ when 
\begin{align}\label{mesh-h-1}
h<(C_2M)^{-4/5} 
\le (C_2\|\u^0\|_{H^2})^{-4/5} .
\end{align}

The induction assumption \eqref{mathind} implies 
$\|\u\|_{L^\infty(0,t_{k};L^4)}
\le \|\u_{h,\tau}\|_{L^\infty(0,t_{k};L^4)}+1 $ and therefore,  
$$
\|\u\|_{L^\infty(0,t_{k};L^4)}
+\|\u^0\|_{H^2} \le M.
$$
When the stepsize $\tau$ satisfies 
\begin{align}\label{meshcd1}
\tau<\varphi(M) , 
\end{align}
the induction assumption \eqref{mathind} and Lemma \ref{prop1} together imply that the strong solution $\u\in L^\infty(0,t_k;\dot{\bf H}^2)\cap W^{1,\infty}(0,t_k;\dot{\bf L}^2)$ can be extended to time $t_k+\tau = t_{k+1}$, i.e.,  
\begin{align}\label{mathind-c1}
\begin{aligned}
&\mbox{problem \eqref{NSPDE1}--\eqref{NSPDE2} with initial value $\u^0$  has a unique strong solution} \\
&\mbox{$\u\in L^\infty(0,t_{k+1};\dot{\bf H}^2)\cap W^{1,\infty}(0,t_{k+1};\dot{\bf L}^2)$,}
\end{aligned}
\end{align}
satisfying the following quantitative estimate:
\begin{align}
&\|\partial_{tt}\u\|_{L^2(0,t_{k+1} ;\dot H^{-1})} 
+\|\partial_{t}\u\|_{L^\infty(0,t_{k+1};L^2)} 
+\|\partial_{t}\u\|_{L^2(0,t_{k+1};H^1)} 
 \nn\\
&\quad 
+ \|\u\|_{L^2(0,t_{k+1} ;H^2)}
+ \|\u\|_{L^\infty(0,t_{k+1} ;H^2)}
+ \|p\|_{L^\infty(0,t_{k+1} ;H^1)}
\nn\\
&\le \Phi(\|\u\|_{L^\infty(0,t_k;L^4)}
+\|\u^0\|_{H^2}) \nn\\
&\le  \Phi(M)   .
\label{LocEst2}
\end{align}
Under this regularity, the solution $\u$ satisfies the variational equations
\begin{align}
& \big(D_{\tau}\u^{n+1} ,{\bf v}_h\big)
-\big(\u^{n+1} ,\u^{n} \cdot \nabla {\bf v}_h\big) \nn\\
&\quad 
+\mu \big(\nabla \u^{n+1} ,\nabla {\bf v}_h\big) 
-\big(p^{n+1} ,\nabla\cdot{\bf v}_h\big)
=\big({\bf E}^{n+1} ,{\bf v}_h\big) 
+\big({\bf F}^{n+1} , \nabla {\bf v}_h\big) ,
\label{ExEq1}\\[5pt]
&\big(\nabla\cdot \u^{n+1},q_h\big)=0,
\label{ExEq2} 
\end{align} 
for $n=0,1,\cdots,k$, where 
\begin{align}
&{\bf E}^{n+1}=D_{\tau}\u^{n+1}-\partial_t\u^{n+1} 
,\\
&{\bf F}^{n+1}= \u^{n+1}\otimes(\u^{n+1}-\u^{n}) ,
\end{align}
are the truncation errors of temporal discretization, satisfying
\begin{align}\label{Estimates-E-F}
&\sum_{n=0}^{k}\tau \|{\bf E}^{n+1}\|_{\dot H^{-1}}^2 
+\sum_{n=0}^{k}\tau \|{\bf F}^{n+1}\|_{L^2}^2  \nn\\
&\le 
C\|D_{\tau}\u^{n+1}-\partial_t\u^{n+1}\|_{L^2(0,t_{k+1};\dot H^{-1})}^2 
+C\|(\u^{n}-\u^{n+1})\otimes \u^{n+1} \|_{L^2(0,t_{k+1};L^2)}^2 
 \nn\\
&\le C\tau^2  \|\partial_{tt}\u \|_{L^2(0,t_{k+1};\dot H^{-1})}^2 
 +C\tau^2 \|\partial_t\u\|_{L^2(0,t_{k+1};L^4)}^2 \|\u\|_{L^\infty(0,t_{k+1};L^4)}^2  \nn\\ 
&\le C(\Phi(M)^2 +\Phi(M)^4)\tau^2  
\nn\\
&\le C\Phi(M)^4\tau^2  
.
\end{align}
where the last inequality uses the property $\Phi(M)\ge 1$  (to simplify the expression). 

Let ${\bf e}_h^{n+1}:=\u_h^{n+1}-\u^{n+1}_{*,h}$ and $\eta_h^{n+1}:=p_h^{n+1}-p^{n+1}_{*,h}$, where $\u^{n+1}_{*,h}$ and $p^{n+1}_{*,h}$ are the Stokes--Ritz projection defined in \eqref{Ritz-simple}. 
The difference between \eqref{FEMEq1}--\eqref{FEMEq2} and \eqref{ExEq1}--\eqref{ExEq2} gives the following two error equations: 
\begin{align}
& \big(D_{\tau}{\bf e}^{n+1}_h ,{\bf v}_h\big)
-\big(\e^{n+1}_h,\u^{n}_h\cdot \nabla {\bf v}_h\big) 
+\mu \big(\nabla \e^{n+1}_h,\nabla {\bf v}_h\big) 
-\big(\eta^{n+1}_h,\nabla\cdot{\bf v}_h\big) \nn\\
&=
-\big({\bf E}^{n+1} ,{\bf v}_h\big) 
-\big({\bf F}^{n+1} ,\nabla {\bf v}_h\big) 
\nn\\
&\quad 
+\big(D_{\tau} \u^{n+1}-D_{\tau} \u^{n+1}_{*,h}) ,{\bf v}_h\big)
\nn \\
&\quad 
-\big(\u^{n+1}-\u^{n+1}_{*,h},\u^{n}_h\cdot \nabla {\bf v}_h\big) \nn\\
&\quad
+\big(\u^{n+1},\e^{n}_h\cdot \nabla {\bf v}_h\big) \nn\\
&\quad 
+\big(\u^{n+1},(\u^{n}_{*,h}-\u^{n})\cdot \nabla {\bf v}_h\big) \nn\\
&=: \sum_{j=1}^5 I_j(\v_h) \qquad\qquad\forall\, \v_h\in {\bf X}_h,
\label{ErrEq1} \\[10pt]
&\big(\nabla\cdot \e^{n+1}_h,q_h\big)=0 \qquad\,\forall\, q_h\in V_h .
\label{ErrEq2} 
\end{align}
Substituting ${\bf v}_h=\e_h^{n+1}$ into \eqref{ErrEq1} and using \eqref{ErrEq2} with $q_h=\eta_h^{n+1}$, we obtain 
\begin{align}\label{energy-error-eq}
D_{\tau}\bigg(\frac{1}{2}\|{\bf e}^{n+1}_h\|_{L^2}^2 \bigg)
+\mu\|\nabla \e^{n+1}_h\|_{L^2}^2 
&=\sum_{j=1}^5 I_j({\bf e}^{n+1}_h) ,
\end{align}
where 
\begin{align*}
I_1({\bf e}^{n+1}_h)
&=
-\big({\bf E}^{n+1} ,{\bf e}_h^{n+1}\big) 
-\big({\bf F}^{n+1} ,\nabla {\bf e}_h^{n+1}\big) \\
&\le
C\big( \|{\bf E}^{n+1}\|_{\dot H^{-1}}
+ \|{\bf F}^{n+1}\|_{L^2} \big)
\|\nabla {\bf e}_h^{n+1}\|_{L^2} \\
&\le
C\mu^{-1}(\|{\bf E}^{n+1}\|_{\dot H^{-1}}^2
+ \|{\bf F}^{n+1}\|_{L^2}^2 ) 
+\frac{\mu}{16}\|\nabla {\bf e}_h^{n+1}\|_{L^2}^2 , \\[10pt]
I_2({\bf e}^{n+1}_h)
&=\big(D_{\tau} \u^{n+1}-R_hD_{\tau} \u^{n+1}) ,{\bf e}^{n+1}_h\big) \\
&\le
Ch\|D_{\tau} \u^{n+1}\|_{H^1} \|{\bf e}^{n+1}_h\|_{L^2} \\
&\le
Ch\|D_{\tau} \u^{n+1}\|_{H^1} \|\nabla {\bf e}^{n+1}_h\|_{L^2} \\
&\le
C\mu^{-1}h^2\|D_{\tau} \u^{n+1}\|_{H^1}^2 
+\frac{\mu}{16} \|\nabla {\bf e}^{n+1}_h\|_{L^2}^2 , \\[10pt]
I_3({\bf e}^{n+1}_h)
&= -\big(\u^{n+1}-\u^{n+1}_{*,h},
\u^{n}_h\cdot \nabla {\bf e}_h^{n+1}\big) \\
&\le 
\|\u^{n+1}-\u^{n+1}_{*,h}\|_{L^6}\|\u_h^{n}\|_{L^3}
\|\nabla {\bf e}_h^{n+1}\|_{L^2} \\
&\le 
C\|\u^{n+1}-\u^{n+1}_{*,h}\|_{H^1}\|\u_h^{n}\|_{L^3}
\|\nabla {\bf e}_h^{n+1}\|_{L^2} \\
&\le 
Ch\|\u^{n+1}\|_{H^2} \|\u_h^{n}\|_{L^3}
\|\nabla {\bf e}_h^{n+1}\|_{L^2} \\
&\le 
C\mu^{-1}h^2 \|\u^{n+1}\|_{H^2}^2 \|\u_h^{n}\|_{L^3}^2 
 +\frac{\mu}{16}
\|\nabla {\bf e}_h^{n+1}\|_{L^2}^2 \\
&\le 
C\mu^{-1}h^2 \|\u^{n+1}\|_{H^2}^2 \|\nabla\u_h^{n}\|_{L^2}^2  
+\frac{\mu}{16}
\|\nabla {\bf e}_h^{n+1}\|_{L^2}^2 , \\[10pt]
I_4({\bf e}^{n+1}_h)
&=\big(\u^{n+1},\e^{n}_h\cdot \nabla {\bf e}^{n+1}_h\big) \\ 
&\le
\|\u^{n+1}\|_{L^6}\|\e^{n}_h\|_{L^3}
\| \nabla {\bf e}^{n+1}_h\|_{L^2} \\
&\le
C\|\nabla \u^{n+1}\|_{L^2}
\|{\bf e}^{n}_h\|_{L^2}^{\frac12}\|{\bf e}^{n}_h\|_{L^6}^{\frac12}
\| \nabla {\bf e}^{n+1}_h\|_{L^2} \\
&\le 
C\|\nabla \u^{n+1}\|_{L^2}
\|{\bf e}^{n}_h\|_{L^2}^{\frac12}\|\nabla {\bf e}^{n}_h\|_{L^2}^{\frac12}
\| \nabla {\bf e}^{n+1}_h\|_{L^2} \\
&\le
C\mu^{-1}\|\nabla \u^{n+1}\|_{L^2}^4 
\|{\bf e}^{n}_h\|_{L^2}^2 
+\frac{\mu}{16} \|\nabla {\bf e}^{n}_h\|_{L^2}^2
+\frac{\mu}{16} \|\nabla {\bf e}^{n+1}_h\|_{L^2}^2, \\[10pt]
I_5({\bf e}^{n+1}_h)
&= \big(\u^{n+1},(\u^n_{*,h}-\u^{n})\cdot \nabla {\bf e}_h^{n+1}\big) \\
&\le 
 \|\u^{n+1}\|_{L^6} \|\u^n_{*,h}-\u^{n}\|_{L^3} \|\nabla {\bf e}_h^{n+1}\|_{L^2} \\
&\le 
C\|\nabla \u^{n+1}\|_{L^2}
\|\u^n_{*,h}-\u^{n}\|_{H^1}
\|\nabla {\bf e}_h^{n+1}\|_{L^2} \\
&\le 
C\|\nabla \u^{n+1}\|_{L^2} 
h \|\u^{n}\|_{H^2} 
\|\nabla {\bf e}_h^{n+1}\|_{L^2} \\
&\le 
C\mu^{-1}h^2 \|\nabla \u^{n+1}\|_{L^2} ^2 \|\u^{n}\|_{H^2} ^2 
+\frac{\mu}{16}
\|\nabla {\bf e}_h^{n+1}\|_{L^2}^2 .
\end{align*}

Substituting the estimates of $I_j({\bf e}^{n+1}_h)$, $j=1,\dots,5$, into \eqref{energy-error-eq}, we obtain
\begin{align}
&\quad D_{\tau}\bigg(\frac{1}{2}\|{\bf e}^{n+1}_h\|_{L^2}^2 \bigg)
+\mu\|\nabla \e^{n+1}_h\|_{L^2}^2 \nn\\
&\quad \le   C_2\mu^{-1} (\|{\bf E}^{n+1}\|_{\dot H^{-1}}^2
+\|{\bf F}^{n+1}\|_{L^2}^2) \nn\\
&\quad \quad 
+C_2\mu^{-1} h^2( \|D_{\tau} \u^{n+1}\|_{H^1}^2+ \|\u^{n+1}\|_{H^2}^2\| \nabla\u_h^{n}\|_{L^2}^2 + \|\u^{n+1}\|_{H^2}^2 \|\nabla \u^{n+1}\|_{L^2} ^2 ) \nn\\
&\quad \quad 
+ C_2\mu^{-1} \|\nabla\u^{n+1}\|_{L^2}^4\| \e_h^{n}\|_{L^2}^2 
+\frac{5\mu}{16}\|\nabla \e_h^{n+1}\|_{L^2}^2 
+\frac{\mu}{16}\|\nabla \e_h^{n}\|_{L^2}^2 . \nn\\[-20pt]
\label{longEst10}
\end{align}\vspace{-15pt}

By using \eqref{LocEst2} and \eqref{Estimates-E-F}, we have 
\begin{align*}
&\sum_{n=0}^{k}\tau \|D_\tau\u^{n+1}\|_{H^1}^2 
\le C\|\partial_t\u\|_{L^2(0,t_{k+1};H^1)}^2 
\le C\Phi(M)^2  , \\[5pt]
& \sum_{n=0}^{k}\tau (\|{\bf E}^{n+1}\|_{\dot H^{-1}}^2 
+ \|{\bf F}^{n+1}\|_{L^2}^2 )
 \le C\Phi(M)^4\tau^2 ,\\[5pt]
& \sum_{n=0}^{k}\tau (  \|\u^{n+1}\|_{H^2}^2\| \nabla\u_h^{n}\|_{L^2}^2 + \|\u^{n+1}\|_{H^2}^2 \|\nabla \u^{n+1}\|_{L^2} ^2 ) \\
&\le \bigg[\sum_{n=0}^{k}\tau ( \| \nabla\u_h^{n}\|_{L^2}^2 +\|\nabla \u^{n+1}\|_{L^2}^2) \bigg]
\|\u \|_{L^\infty(0,t_{k+1};H^2)}^2 \\
&\le \bigg[C\| \u_h^0\|_{L^2}^2+C\sum_{n=0}^{k}\tau  \|\nabla \u^{n+1}\|_{L^2}^2   \bigg]
\Phi(M)^2  \qquad\mbox{(here \eqref{EnergEst} is used)} \\
&\le \bigg[C\| \u^{0}\|_{H^2}^2
+C\|\nabla \u \|_{L^2(0,t_{k+1};L^2)}^2 
+C\tau^2 \|\partial_t\nabla\u\|_{L^2(0,t_{k+1};L^2)}^2   \bigg]\Phi(M)^2  \\
&\le C\Phi(M)^4 ,
\end{align*}
where we have used the expression 
$$\nabla\u^{n+1}=\frac{1}{\tau}\int_{t_n}^{t_{n+1}}\nabla\u(t)\d t
+\frac{1}{\tau}\int_{t_n}^{t_{n+1}}(s-t_n)\partial_t\nabla\u(s)\d s$$
in estimating
$\sum_{n=0}^{k}\tau  \|\nabla \u^{n+1}\|_{L^2}^2  $.
By using the last three estimates
and summing up \eqref{longEst10} for $n=0,1,\cdots,m$ 
(with $0\le m\le k$), we obtain
\begin{align}\label{longEst1}
\frac{1}{2}\|{\bf e}^{m+1}_h\|_{L^2}^2  
+\frac{\mu}{2}\sum_{n=0}^m\tau
\|\nabla \e^{n+1}_h\|_{L^2}^2 
&\le \frac{1}{2}\|{\bf e}^0_h\|_{L^2}^2  
+\frac{\tau\mu}{16}\|\nabla{\bf e}^0_h\|_{L^2}^2
+C\mu^{-1} \Phi(M)^4(\tau^2+h^2) \nn\\
&\quad
+ C\mu^{-1} \sum_{n=0}^m\tau \|\nabla \u^{n+1}\|_{L^2}^4 \|\e_h^{n}\|_{L^2}^2 . 
\end{align}

Since $\u_h^0$ is the Stokes--Ritz projection of $\u^0$, it follows that the initial error satisfies 
\begin{align*}
\|\e_h^0\|_{L^2} + \tau\mu \|\nabla \e_h^0\|_{L^2} 
\le 
C\|\u^0\|_{H^2} h +  C\|\u^0\|_{H^2} \tau 
&\le C\Phi(M) (\tau+h) \\ 
&\le C\Phi(M)^2 (\tau+h) ,
\end{align*}
where we have used $\Phi(M)\ge 1$ in the last inequality. Furthermore, \eqref{LocEst2} implies  
$$
\|\nabla \u^{n+1}\|_{L^2}^2
\le
\|\u\|_{L^\infty(0,t_{k+1};H^1)}^2
\le 
\Phi(M)^2 \quad\mbox{for}\,\,\,  0\le n\le k . 
$$ 
By using the last two estimates, \eqref{longEst1} can be reduced to 
\begin{align}
&\frac{1}{2}\|{\bf e}^{m+1}_h\|_{L^2}^2  
+\frac{\mu}{2}\sum_{n=0}^m\tau
\|\nabla \e^{n+1}_h\|_{L^2}^2 \nn\\
&\le C\mu^{-1}  \Phi(M)^4(\tau^2+h^2) 
+C_3\mu^{-1} \Phi(M)^2\sum_{n=0}^m\tau \|\nabla \u^{n+1}\|_{L^2}^2 (\|\e_h^{n+1}\|_{L^2}^2 + \|\e_h^{n}\|_{L^2}^2) .
\label{longEst}
\end{align}
Then, using the expression 
$\nabla\u^{n+1}=\frac{1}{\tau}\int_{t_n}^{t_{n+1}}\nabla\u(t)\d t
+\frac{1}{\tau}\int_{t_n}^{t_{n+1}}(s-t_n)\partial_t\nabla\u(s)\d s$, we have 
\begin{align}
\sum_{n=0}^{k }\tau \|\nabla\u^{n+1}\|_{L^2}^2 
&\le 2\|\nabla\u \|_{L^2(0,t_{k+1};L^4)}^2 
+2\tau^2 \|\partial_t\nabla\u\|_{L^2(0,t_{k+1};L^4)}^2 \nn\\
&\le C\|\nabla\u \|_{L^2(0,t_{k+1};H^1)}^2 
+C\tau^2 \|\partial_t\nabla\u\|_{L^2(0,t_{k+1};H^1)}^2 \nn\\
&\le C\Phi(M)^2 ,
\end{align}
where the last inequality is again due to \eqref{LocEst2}. 
Applying Gronwall's inequality to \eqref{longEst} and using the inequality above, we obtain 
\begin{align}
\max_{0\le n\le k}\|{\bf e}^{n+1}_h\|_{L^2}^2 
&\le  \exp\left(C\mu^{-1} \Phi(M)^2
\sum_{n=0}^{k}\tau \|\nabla \u^{n+1}\|_{L^2}^2 \right) C\mu^{-1}  \Phi(M)^4 (\tau^2+h^2) \nn\\
&\le \exp(C\mu^{-1} \Phi(M)^4)C\mu^{-1}  \Phi(M)^4 (\tau^2+h^2)  \nn\\
&\le \exp(C_4\mu^{-1} \Phi(M)^4) (\tau^2+h^2) . 
\end{align}
Then, substituting the estimate above into \eqref{longEst}, we obtain the following error estimate:
\begin{align}\label{L2H1Err}
&\max_{0\le n\le k}\|{\bf e}^{n+1}_h\|_{L^2}^2 +   \sum_{n=0}^{k}\tau\|\nabla \e^{n+1}_h\|_{L^2}^2
\le  \exp(C_5\mu^{-1}\Phi(M)^4)(\tau^2+h^2 ) .
\end{align}

By using the inverse inequality 
$\|{\bf e}_h^{n+1}\|_{L^4}
\le
C h^{-3/4}\|{\bf e}_h^{n+1}\|_{L^2} $ and the Sobolev embedding inequality $\|{\bf e}_h^{n+1}\|_{L^4}
\le
C\|\nabla {\bf e}_h^{n+1}\|_{L^2} $, we have 
$$
\|{\bf e}_h^{n+1}\|_{L^4} 
\le 
\min(C h^{-3/4}\|{\bf e}_h^{n+1}\|_{L^2},C\|\nabla {\bf e}_h^{n+1}\|_{L^2}). 
$$
Let ${\bf e}_{h,\tau}(x,t)={\bf e}_h^n(x)$ for $t\in(t_{n-1},t_n]$, as defined in \eqref{DefNSol}. Then, by using the inequality above, from \eqref{L2H1Err} we derive that 
\begin{align}
\max_{0\le n\le k}\|\u_h^{n+1}-\u^{n+1}_{*,h}\|_{L^4}^2  
&= 
\max_{0\le n\le k}\|{\bf e}_{h}^{n+1}\|_{L^4}^2 \nn\\
&\le \min(C_6 h^{-3/2}\max_{0\le n\le k}\|{\bf e}_{h}^{n+1}\|_{L^2}^2,
C_6 \max_{0\le n\le k}\|\nabla \e_{h}^{n+1}\|_{L^2}^2) \nn\\
&\le C_6 \min( h^{-3/2}  ,\tau^{-1}) \biggl(\max_{0\le n\le k}\|{\bf e}^{n+1}_h\|_{L^2}^2 + \sum_{n=0}^{k}\tau\|\nabla \e^{n+1}_h\|_{L^2}^2\bigg)\nn\\
&\le C_6\exp(C_5\mu^{-1}\Phi(M)^4)(\tau+h^{1/2} ) .
\end{align}
For any $t\in (t_n,t_{n+1}]$ and $n=0,\dots,k$, we have
\begin{align}
&\max_{t\in (t_n,t_{n+1}]}
\|\u^{n+1}_{*,h}-\u(\cdot,t)\|_{L^4}^2 \nn\\ 
&\le 2\|\u^{n+1}_{*,h}-\u^{n+1}\|_{L^4}^2
+\max_{t\in (t_n,t_{n+1}]} 2\|\u^{n+1}-\u(\cdot,t)\|_{L^4}^2 \nn \\
&\le 2\|\u^{n+1}_{*,h}-\u^{n+1}\|_{H^1}^2
+\max_{t\in (t_n,t_{n+1}]} 2\|\u^{n+1}-\u(\cdot,t)\|_{L^4}^2 \nn \\
&\le C \|\u^{n+1}\|_{H^2}^2 h^2 
+ 2 \tau^2\|\partial_t\u\|_{L^\infty(t_n,t_{n+1};L^4)}^2\nn\\
&\le C \|\u\|_{L^\infty(0,t_{k+1};H^2)}^2 h^2 
+ C \tau^2\|\partial_t\u\|_{L^\infty(0,t_{k+1};H^1)}^2\nn\\
&\le C\Phi(M)^2 (\tau^2+h^2) \nn\\
&\le C\mu^{-1}\Phi(M)^4 (\tau^2+h^2) \nn\\
&\le \exp(C_7\mu^{-1}\Phi(M)^4)(\tau^2+h^2) .
\end{align}
Combining the two estimates above and using the triangle inequality, we obtain 
\begin{align}
&\|\u_{h,\tau}-\u\|_{L^\infty(0,t_{k+1};L^4)}^2\leq
\exp(C_8\mu^{-1}\Phi(M)^4)(\tau+h^{1/2})  
\end{align} 
for some positive constant $C_8$.
When 
\begin{align}\label{meshcd3}
\tau+h^{1/2} \le \exp(-C_8\mu^{-1}\Phi(M)^4),
\end{align}
we have
\begin{align}\label{mathind-c2}
\|\u_{h,\tau}-\u\|_{L^\infty(0,t_{k+1};L^4)}  \le 1 .
\end{align}
Since the right-hand side of \eqref{meshcd3} is independent of $k$ (depending only on $M$), the existence of strong solution in \eqref{mathind-c1} and estimate \eqref{mathind-c2} together complete the mathematical induction on \eqref{mathind}.

Overall, if the mesh conditions \eqref{mesh-h-1}--\eqref{meshcd1} and \eqref{meshcd3} are satisfied, then by mathematical induction on $k$ in \eqref{mathind} for $k=1,\dots,N$, we have proved the existence of a unique strong solution $\u\in L^\infty(0,T;\dot{\bf H}^2)\cap W^{1,\infty}(0,T;\dot{\bf L}^2)$, satisfying 
\begin{align}
\|\u\|_{L^\infty(0,T;L^4)}  \le \|\u_{h,\tau}\|_{L^\infty(0,T;L^4)}+1   .
\end{align}
Therefore, Theorem \ref{MainTHM} is proved with
\begin{align}\label{expr-tauM}
\tau_M=\min\bigg(\varphi(M) ,  \frac{1}{2}\exp(-C_8\mu^{-1}\Phi(M)^4)\bigg)
\end{align}
and
\begin{align}\label{expr-hM}
h_M= \min\bigg(\frac{1}{(C_2M)^{\frac45}},\frac{1}{2}\exp(-C_8\mu^{-1}\Phi(M)^4) \bigg) .
\end{align}

\begin{remark}$\,$
Clearly, we can choose $C_8\geq C_5$ in the analysis above. In this case, \eqref{L2H1Err} and \eqref{meshcd3} imply an error estimate:
\begin{align}\label{L2H1ErrF}
&\|\u_{h,\tau}-\u\|_{L^\infty(0,T;L^2)}^2
\le 2(\tau +h^{3/2}) .
\end{align}
\end{remark}


\bigskip
 
\begin{thebibliography}{99}

\bibitem{Adams}
R. A. Adams and J. J. F. Fournier, {\em Sobolev spaces}, Second Edition, Academic Press, 
The Netherlands, 2003.

\bibitem{B73}
I. Babuska, {\em The finite element method with Lagrangian multipliers}, Numer. Math. 20 (1973), pp. 179--192.

\bibitem{BSS12}
R. Bermejo, P. G. del Sastre, and L. Saavedra,
{\em A second order in time modified Lagrange--Galerkin finite element method for the incompressible Navier--Stokes equations}, SIAM J. Numer. Anal. 50 (2012), pp. 3084--3109.

\bibitem{Braess}
D. Braess, {\em Finite Elements: Theory, Fast Solvers, and Applications in Solid Mechanics,} 
Cambridge University Press, Cambridge, UK, 1997.

\bibitem{Buckmaster-Vicol-2019}
T. Buckmaster and V. Vicol, 
{\em Nonuniqueness of Weak Solutions to the Navier--Stokes Equation}, 
Annals of Mathematics 189 (2019), pp. 101--144. 

\bibitem{Brenner-Scott}
S. C. Brenner and L. R. Scott,   
{\em The Mathematical Theory of Finite Element Methods},
Third edition. Texts in Applied Mathematics 15, Springer, New York, 2008. 

\bibitem{B74}
F. Brezzi, {\em On the existence, uniqueness and approximation of saddle point problem sarising from Lagrangian multipliers}, RAIRO. Anal. Numer. 8 (1974), pp. 129--151.

%

\bibitem{Constantin-Foias}
P. Constantin and C. Foias. 
{\em Navier--Stokes Equations}, The University of Chicago Press, Chicago, 1988.

\bibitem{Dauge-1989}
M. Dauge, 
{\em Stationary Stokes and Navier--Stokes systems on two- or three-dimensional domains with corners. Part I. Linearized equations}, SIAM J. Math. Anal. 20 (1989), pp. 74--97.

\bibitem{FK64}
H. Fujita and T. Kato, 
{\em On the Navier-Stokes initial value problem I}, 
Arch. Rational Mech. Anal. 16 (1964), pp. 269--315.

\bibitem{Cheskidov}
A. Cheskidov, 
{\em Blow-up in finite time for the dyadic model of the 
Navier--Stokes equations}, 
Trans. Am. Math. Soc. 360 (2010), pp. 5101--5120.

\bibitem{Chorin68}
A. J. Chorin,
{\em Numerical solution of the Navier--Stokes equations},
Math. Comp. 22 (1968), pp. 745--762.

\bibitem{Chorin}
A. J. Chorin, 
{\em On the convergence of discrete approximations to the Navier--Stokes equations}, 
Math. Comp. 23 (1969), pp. 341--353.

\bibitem{DGLQ82}
J. Donea, S. Giuliani, H. Laval, and L. Quartapelle,
{\em Finite element solution of the unsteady Navier--Stokes equations by a fractional step method}, 
Comput. Meths. Appl. Mech. Eng. 30 (1982), pp. 53--73.

\bibitem{GP09}
I. Gallagher and M. Paicu, 
{\em Remarks on the blow-up of solutions to a toy model for the Navier-Stokes equations}, 
Proc. Amer. Math. Soc. 137 (2009), pp. 2075--2083.

\bibitem{Giga}
Y. Giga, 
{\em Solutions for semilinear parabolic equations in Lp and regularity of weak solutions of the Navier--Stokes system}, 
J. Differ. Equations 62 (1986), pp. 186--212. 

\bibitem{2014-Guzman-Neilan}
J. Guzm\'an and M. Neilan, 
{\em Conforming and divergence-free Stokes elements on general triangular
  meshes}, Math. Comp. 83 (2014), pp. 15--36.

\bibitem{2019-Guzman-Scott}
J. Guzman and L. R. Scott, 
{\em The Scott-Vogelius finite elements revisited}, 
Math. Comp. 88 (2019), pp. 515--529. 

\bibitem{Heywood82}
J. G. Heywood, 
{\em An error estimate uniform in time for 
spectral Galerkin approximations of the Navier--Stokes problem},
Pacific J. Math. 98 (1982), pp. 333--345.

\bibitem{HR82}
J. G. Heywood and R. Rannacher, 
{\em Finite element approximation of the nonstationary Navier--Stokes problem. I. Regularity of solutions and second-order spatial discretization}, 
SIAM J. Numer .Anal. 19 (1982), pp. 275--311.

\bibitem{HR90}
J. G. Heywood and  R. Rannacher, 
{\em Finite-element approximation of the nonstationary 
Navier--Stokes problem Part IV: error analysis for second-order time discretization}, 
SIAM J. Numer. Anal. 27 (1990), pp. 353--384.

\bibitem{Hopf51}
E. Hopf, 
{\em \"Uber die Anfangswertaufgabe f\"ur die hydrodynamischen Grundgleichungen}, 
Math. Nachr. 4 (1951), pp. 213--231.

\bibitem{HL09}
T. Hou and Z. Lei, 
{\em On the stabilizing effect of convection in three-dimensional incompressible flows},
Comm. Pure Appl. Math., 62 (2009), pp. 501--564.

\bibitem{He-2013}
Y. He, 
{\em Euler implicit/explicit iterative scheme for the stationary Navier--Stokes equations}, 
Numer. Math. 123 (2013), pp. 67--96.

\bibitem{Kato84}
T. Kato, 
{\em Strong $L^p$ solutions of the Navier-Stokes equations in $\R^m$, with applications to weak solutions}. 
Math. Z. 187 (1984), pp. 471--480.


\bibitem{KT01}
H. Koch and D. Tataru,
{\em Well-posedness for the Navier-Stokes equations},
Adv. Math. 157 (2001), pp. 22--35.
 
\bibitem{LL11}
Z. Lei and F. Lin,
{\em Global mild solutions of Navier-Stokes equations},
Comm. Pure Appl. Math. LXIV (2011),  pp. 1297--1304.

\bibitem{Leray34}
J. Leray, {\em Sur le mouvement d'un liquide visqueux emplissant l'espace}, 
Acta Math. 63 (1934), pp. 193--248.

\bibitem{Lions}
P. L. Lions, {\em Mathematical Topics in Fluid Mechanics,
Volume 1: Incompressible Models}, 
Clarendon Press, Oxford, 1996.

\bibitem{MS01}
S. Montgomery-Smith, 
{\em Finite time blow up for a Navier--Stokes like equation}, 
Proc. Amer. Math. Soc. 129 (2001), pp. 3025--3029.

\bibitem{Pironneau}
O. Pironneau, 
{\em On the transport-diffusion 
algorithm and its applications to the Navier--Stokes
equations}, 
Numer. Math. 38 (1982), pp. 309--332.

\bibitem{Rann00}
R. Rannacher,
{\em Finite element methods for the incompressible 
Navier--Stokes equations}, 
in ``Fundamental Directions in Mathematical Fluid Mechanics'', pp. 191--293, edited by G.P. Galdi, J. Heywood and R. Ranancher. Birkh\"auser, Basel, Boston, Berlin, 2000.

\bibitem{Serrin}
J. Serrin, 
{\em On the interior regularity of weak solutions of the Navier--Stokes equations}, Arch. Rational Mech. Anal. 9 (1962), pp. 187--195.

\bibitem{Shen92}
J. Shen,
{\em On error estimates of the projection methods for the Navier-Stokes equations: first-order schemes}, 
SIAM J. Numer. Anal. 29 (1992), pp. 57--77. 

\bibitem{STW}
J. Shen, T. Tang, and L. Wang,
{\em Spectral Methods: Algorithms, Analysis and Applications},
Springer-Verlag, Berlin-Heidelberg, 2011.


\bibitem{Stein1970}
E. M. Stein, 
{\em Singular Integrals and Differentiability Properties of Functions}, 
Princeton University Press, New Jersey, 1970. 

\bibitem{Suli88}
E. S\"uli, 
{\em Convergence and non-linear stability of the Lagrange--Galerkin method for the Navier--Stokes equations}, 
Numer. Math. 53 (1988), pp. 459--483.

\bibitem{Tao}
T. Tao, 
{\em Finite time blowup for an averaged three-dimensional Navier--Stokes equation}, 
J. Amer. Math. Soc. 29 (2016), pp. 601--674. 

\bibitem{Taylor}
M. E. Taylor, 
{\em Incompressible fluid flows on rough domains},  
Progress in Nonlinear Differential Equations and Their Applications, 42 (2000), pp. 320--334.

\bibitem{Temam}
R. Temam, 
{\em Navier-Stokes Equations, 
Theory and Numerical Analysis}, 1979, North-Holland.

\bibitem{Triebel}
H. Triebel, 
{\em Interpolation Theory, Function Spaces, Differential Operators}, 
Amsterdam-New York-Oxford: North-Holland, 1978. 

\bibitem{Verfurth}
R. Verf\"urth,
{\em Error estimates for a mixed finite element 
approximation of the Stokes equations}, 
RAIRO. Anal. Numer. 18 (1984), pp. 175--182.

\end{thebibliography}
\end{document}